# EVOLUTIONARILY STABLE STRATEGIES OF RANDOM GAMES, AND THE VERTICES OF RANDOM POLYGONS[1]

By Sergiu Hart, Yosef Rinott and Benjamin Weiss

*Hebrew University of Jerusalem*

An *evolutionarily stable strategy* (ESS) is an equilibrium strategy that is immune to invasions by rare alternative ("mutant") strategies. Unlike Nash equilibria, ESS do not always exist in finite games. In this paper we address the question of what happens when the size of the game increases: does an ESS exist for "almost every large" game? Letting the entries in the $n \times n$ game matrix be independently randomly chosen according to a distribution $F$, we study the number of ESS with support of size 2. In particular, we show that, as $n \to \infty$, the probability of having such an ESS: (i) converges to 1 for distributions $F$ with "exponential and faster decreasing tails" (e.g., uniform, normal, exponential); and (ii) converges to $1 - 1/\sqrt{e}$ for distributions $F$ with "slower than exponential decreasing tails" (e.g., lognormal, Pareto, Cauchy).

Our results also imply that the expected number of vertices of the convex hull of $n$ random points in the plane converges to infinity for the distributions in (i), and to 4 for the distributions in (ii).

**1. Introduction.** The concept of *evolutionarily stable strategy* (ESS for short), introduced by Maynard Smith and Price [12], refers to a strategy that, when played by the whole population, is immune to invasions by rare alternative ("mutant") strategies (see Section 2.1 for precise definitions). Formally, an ESS corresponds to a symmetric Nash equilibrium that satisfies an additional stability requirement. Every (symmetric) finite game has a (symmetric) Nash equilibrium. But the same is not true for ESS: there are games with finitely many pure strategies that have no ESS. Moreover, the

Received February 2007; revised June 2007.

[1]Supported in part by grants of the Israel Science Foundation and by the Institute for Advanced Studies at the Hebrew University of Jerusalem.

nonexistence of ESS is not an "isolated" phenomenon: it holds for open sets of games.[2]

This leads us to the question of what happens when the number of strategies is large: does an ESS exist for "almost every large game"? Specifically, assuming that the payoffs in the game are randomly chosen (they are independent and identically distributed random variables), what is the probability that an ESS exists, and what is the limit of this probability as the size of the game increases?

For *pure ESS*, the answer to this question is simple: the probability that a pure ESS exists is $1 - (1 - 1/n)^n$, which converges to $1 - 1/e \simeq 63\%$ as $n \to \infty$, where $n$ is the number of strategies. What about mixed ESS? Here we study mixed ESS with support of size 2—called "*two-point ESS*"—and find out that, unlike pure ESS, the answer depends on the underlying distribution $F$ from which the payoffs are drawn.

By way of illustration, consider the family of cumulative distribution functions $F_\alpha(x) = 1 - e^{-x^\alpha}$ for all $x \geq 0$, where $\alpha > 0$. Our result is:

- When $\alpha \geq 1$ the probability that there is a two-point ESS converges to 1 as[3] $n \to \infty$.
- When $\alpha < 1$ the probability that there is a two-point ESS converges to $1 - 1/\sqrt{e} \simeq 39\%$ as[4] $n \to \infty$.

Moreover, we show that the distribution of the number of two-point ESS converges to a Poisson distribution, with a parameter converging to infinity when $\alpha \geq 1$, and with a parameter of $1/2$ when $\alpha < 1$.

This threshold phenomenon is not restricted to the class $F_\alpha$. We identify two classes of distributions. The first is a class of "light-tailed" distributions with tail probabilities $1 - F(x)$ that decrease exponentially as $x \to \infty$ (i.e., exponential distributions) or faster (e.g., normal distributions, uniform distributions on bounded intervals, logistic distributions); they all lead to the same result as $F_\alpha$ for $\alpha \geq 1$. The second is a class of "heavy-tailed" distributions with tail probabilities that decrease slower than exponentially as $x \to \infty$ (including, in particular, the following distributions: Pareto, Cauchy, lognormal, stable with parameter less than 2), which all behave like $F_\alpha$ for $\alpha < 1$. We refer to these two classes, respectively, as $\mathcal{EF}$ for "*Exponential and Faster* decreasing tails," and $\mathcal{SE}$ for "*Slower than Exponential* decreasing tails" (see Sections 4 and 5 for precise definitions).

---

[2]For instance, the "rock-scissors-paper" game of Example 9.2.1 in van Damme [15], and all its small enough perturbations, have no ESS.

[3]So a fortiori the probability that an ESS exists converges to 1 in this case.

[4]We also show in this case that the probability that there is either a pure or a two-point ESS converges to $1 - e^{-3/2} \simeq 78\%$.

An interesting consequence of our results concerns the classic problem of the number of vertices of the convex hull of a collection of random points in the plane, originally studied by Rényi and Sulanke [13]; see Section 3. Taking symmetric versions of the distributions[5] $F_\alpha$, and assuming that the $2n$ coordinates of the $n$ points in the plane are independent and $F_\alpha$-distributed, we have:

- When $\alpha \geq 1$ the expected number of vertices of the convex hull of $n$ random points in the plane converges to *infinity* as $n \to \infty$.
- When $\alpha < 1$ the expected number of vertices of the convex hull of $n$ random points in the plane converges to 4 as $n \to \infty$.

In addition, in the second case $\alpha < 1$, the number of vertices converges in probability to 4; thus, the convex hull is a quadrilateral with probability converging to 1. Here again, the results hold for the general classes $\mathcal{FE}$ and $\mathcal{SE}$, respectively.

The paper is organized as follows. The two classes of distributions are defined in Sections 4.1 and 5.1, respectively. Our main results for ESS are stated in Theorems 1 and 2 in Section 2.2 (see also Theorem 17 in Section 4.2 and Theorem 33 in Section 5.3), and, for the number of vertices, in Theorem 10 in Section 3. Section 2 presents the model—ESS and random games—together with some preliminary results. Section 3 deals with the number of vertices of random polygons. The detailed analysis is provided in Sections 4 and 5, and we conclude with a discussion in Section 6.

## 2. Preliminaries.

2.1. *Evolutionarily stable strategies.* The setup is that of a symmetric two-person game, with the payoffs given by the $n \times n$ matrix $R = (R(i,j))_{i,j=1,\ldots,n}$. The interpretation is that a meeting between two players, the first playing the pure strategy $i$ and the second playing the pure strategy $j$ (where $1 \leq i, j \leq n$), yields a payoff of $R(i,j)$ to the first, and $R(j,i)$ to the second (these payoffs may be viewed as a measure of "fitness" or "reproductive success").[6] A mixed strategy $p$ is a probability vector on the set of pure strategies, that is, $p = (p_1, \ldots, p_n) \in \Delta(n) := \{x \in \mathbb{R}_+^n : \sum_{i=1}^n x_i = 1\}$; the payoff function $R$ is bilinearly extended to pairs of mixed strategies: $R(p,q) := \sum_{i=1}^n \sum_{j=1}^n p_i q_j R(i,j)$.

A mixed strategy $p \in \Delta(n)$ is an *evolutionarily stable strategy* (*ESS*) for the matrix $R$ if it satisfies the following conditions (Maynard Smith and Price [12]):

---

[5]That is, $F_\alpha(x) = (1/2)e^{-|x|^\alpha}$ for $x \leq 0$ and $F_\alpha(x) = 1 - (1/2)e^{-x^\alpha}$ for $x \geq 0$ [a distribution $F$ is *symmetric* if $F(-x) = 1 - F(x)$ for all $x$].

[6]Thus the payoff matrix of the first player is $R$, and that of the second player is $R^\top$, the transpose of $R$.

[ESS1] $R(p, p) \geq R(q, p)$ for all $q \in \Delta(n)$.

[ESS2] If $q \neq p$ satisfies $R(q, p) = R(p, p)$, then $R(q, q) < R(p, q)$.

This definition is equivalent to the requirement that for every $q \neq p$ there exists an "invasion barrier" $b(q) > 0$ such that $R(p, (1-\varepsilon)p + \varepsilon q) > R(q, (1-\varepsilon)p + \varepsilon q)$ for all $\varepsilon \in (0, b(q))$. The interpretation of this inequality is that any small enough proportion $\varepsilon$ [i.e., less than $b(q)$] of $q$-mutants cannot successfully invade a $p$-population, since the mutants' (average) payoff is strictly less than that of the existing population.

An ESS $p$ is called an *$\ell$-point ESS* if the support $\text{supp}(p) = \{i : p_i > 0\}$ of $p$ is of size $\ell$. In particular, when $\ell = 1$ we have a *pure ESS*. In the biological setup, $\ell = 1$ corresponds to "monomorphism," and $\ell > 1$ to "$\ell$-allele polymorphism." Let $S_\ell^{(n)} \equiv S_\ell^{(n)}(R)$ be the number of $\ell$-point ESS for the matrix $R$.

### 2.2. *ESS of random games.*

Let $F$ be a cumulative distribution function on $\mathbb{R}$. We will assume throughout this paper that $F$ is continuous with a support $(a, b)$ that is either finite or infinite (i.e., $-\infty \leq a < b \leq \infty$). For every integer $n \geq 1$, let $R \equiv R^{(n)}$ be an $n \times n$ matrix whose $n^2$ elements are independent $F$-distributed random variables; the number of $\ell$-point ESS of $R^{(n)}$ is now a random variable $S_\ell^{(n)}$.

We use the following notation: $\mathbf{E}$ for expectation; $\mathcal{L}(Z)$ for the distribution function of the random variable $Z$; $\text{Poisson}(\lambda)$ for the Poisson distribution with parameter $\lambda$ [i.e., $\mathcal{L}(Z) = \text{Poisson}(\lambda)$ if $\mathbf{P}(Z = k) = e^{-\lambda}\lambda^k/k!$ for all integers $k \geq 0$]; and the convergence of distributions is with respect to the variation norm [i.e., the $l_1$-norm on measures: $\|\mathcal{L}(Z_1) - \mathcal{L}(Z_2)\| = \sum_k |\mathbf{P}(Z_1 = k) - \mathbf{P}(Z_2 = k)|$]. The two classes of distributions, namely, the "exponential and faster decreasing tails" class $\mathcal{EF}$ and the "slower than exponential decreasing tails" class $\mathcal{SE}$, will be formally defined in Sections 4.1 and 5.1, respectively.

We now state our main results on $S_2^{(n)}$, the number of two-point ESS:

THEOREM 1. *If $F \in \mathcal{EF}$, then, as $n \to \infty$:*

(i) $\mu_n := \mathbf{E}(S_2^{(n)}) \to \infty$;

(ii) $\|\mathcal{L}(S_2^{(n)}) - \text{Poisson}(\mu_n)\| \to 0$; *and*

(iii) $\mathbf{P}(\text{there is a two-point ESS}) \to 1$.

THEOREM 2. *If $F \in \mathcal{SE}$, then, as $n \to \infty$:*

(i) $\mu_n := \mathbf{E}(S_2^{(n)}) \to 1/2$;

(ii) $\|\mathcal{L}(S_2^{(n)}) - \text{Poisson}(1/2)\| \to 0$; *and*

(iii) $\mathbf{P}(\text{there is a two-point ESS}) \to 1 - e^{-1/2} \simeq 0.39$.

For the convergence to Poisson distributions (ii) we will use a result of the so-called "Chen–Stein method" that requires estimating only the first two moments (see Section 2.5); surprisingly, our proofs in the two cases are different. As for (iii), they are immediate from (ii). The two theorems are proved in Sections 4 and 5, respectively. Note that, for distributions in $\mathcal{EF}$, Theorem 1(iii) implies that the probability that there is an ESS converges to 1 [see Section 6(c)].

Returning to the definition of ESS in Section 2.1, condition [ESS1] says that $p$ is a best reply to itself, that is, $(p, p)$ is a Nash equilibrium. By the bilinearity of $R$, it is equivalent to: $R(i, p) = R(p, p)$ for all $i \in \operatorname{supp}(p)$, and $R(j, p) \leq R(p, p)$ for all $j \notin \operatorname{supp}(p)$. Since $F$ is a continuous distribution, it follows that, with probability 1, the inequalities are strict, that is, $R(j, p) < R(p, p)$ for all $j \notin \operatorname{supp}(p)$ [the $j$th row is independent of the rows in $\operatorname{supp}(p)$]. Therefore, there are no best replies to $p$ outside the support of[7] $p$, that is, $R(q, p) = R(p, p)$ if and only if $\operatorname{supp}(q) \subset \operatorname{supp}(p)$. Thus condition [ESS2] applies only to such $q$, and we obtain (see Haigh [10]):

LEMMA 3. *For a random matrix $R$, the following hold a.s.:*

(i) *$i$ is a pure ESS if and only if $R(i, i) > R(j, i)$ for all $j \neq i$.*

(ii) *There is a two-point ESS with support $\{i, j\}$ if and only if there exist $p_i, p_j > 0$ such that $p_i R(i, i) + p_j R(i, j) = p_i R(j, i) + p_j R(j, j) > p_i R(k, i) + p_j R(k, j)$ for all $k \neq i, j$, and $R(i, i) < R(j, i)$ and $R(j, j) < R(i, j)$.*

The following is immediate from (i) (see Haigh [10]):

PROPOSITION 4. *$S_1^{(n)}$, the number of pure ESS, is a $\operatorname{Binomial}(n, 1/n)$ random variable, and thus $\mathcal{L}(S_1^{(n)}) \to \operatorname{Poisson}(1)$ as $n \to \infty$.*

PROOF. $S_1^{(n)} = \sum_{i=1}^{n} C_i$ where $C_i$ is the indicator that $i$ is a pure ESS, that is, $R(i, i) > R(j, i)$ for all $j \neq i$, and so $\mathbf{P}(C_i = 1) = 1/n$. □

For two-point ESS, we can express their number $S_2^{(n)}$ as a sum of $n(n-1)/2$ identically distributed indicators,

$$S_2^{(n)} = \sum_{1 \leq i < j \leq n} D_{ij},$$

where $D_{ij} \equiv D_{ij}^{(n)}$ is the indicator that columns $i, j$ provide a two-point ESS.[8] To study the asymptotic behavior of $S_2^{(n)}$, we will need to evaluate

---

[7] So $(p, p)$ is a *quasi-strict* Nash equilibrium.

[8] Lemma 3(ii) implies that, a.s., for each $i \neq j$ there can be at most one ESS with support $\{i, j\}$ (in fact, condition [ESS2] implies that the supports of two distinct ESS $p$

the first two moments (see Section 2.5), namely, $\mathbf{P}(D_{ij} = 1) = \mathbf{P}(D_{12} = 1)$ and $\mathbf{P}(D_{ij} = D_{ij'} = 1) = \mathbf{P}(D_{12} = D_{13} = 1)$ (when $\{i, j\}$ and $\{i', j'\}$ are disjoint, $D_{ij}$ and $D_{i'j'}$ are independent, since $D_{ij}$ is a function of the entries in columns $i$ and $j$ only).

2.3. *First moment.* The event that $D_{12} = 1$ depends only on the entries in the first two columns of the matrix $R$, which we will denote $X_i = R(i, 1)$ and $Y_i = R(i, 2)$. Thus $X_1, \ldots, X_n, Y_1, \ldots, Y_n$ are $2n$ independent $F$-distributed random variables. For each $i$, let $P_i := (X_i, Y_i)$ be the corresponding point in $\mathbb{R}^2$. The two points $P_1$ and $P_2$ are almost surely distinct, and thus determine a line $Ax + By = C$ through them, where[9]

$$(1) \qquad A := Y_1 - Y_2, \qquad B := X_2 - X_1, \qquad C := X_2 Y_1 - X_1 Y_2.$$

Finally, we denote by $\Gamma \equiv \Gamma^{(n)}$ the event that there is a two-point ESS with support $\{1, 2\}$, that is, $D_{12} = 1$; recalling Lemma 3(ii), we have

$$\Gamma \equiv \Gamma^{(n)} := \{X_1 < X_2, Y_1 > Y_2, AX_k + BY_k < C \text{ for all } k = 3, \ldots, n\}.$$

Let $\mu_n := \mathbf{E}(S_2^{(n)})$ denote the expected number of two-point ESS. Then

$$(2) \qquad \mu_n = \binom{n}{2} \mathbf{P}(\Gamma^{(n)}).$$

We now define an auxiliary random variable $U \equiv U^{(n)}$, a function of $P_1$ and $P_2$, as follows:

$$(3) \qquad U := \begin{cases} \mathbf{P}(AX_3 + BY_3 > C | P_1, P_2), & \text{if } X_1 < X_2 \text{ and } Y_1 > Y_2, \\ 1, & \text{otherwise,} \end{cases}$$

where $A, B$ and $C$ are determined as above (1) by $P_1$ and $P_2$. Thus $U$ is the probability that an independent point lies above the line through $P_1$ and $P_2$ when $X_1 < X_2$ and $Y_1 > Y_2$. Let $F_U$ be the cumulative distribution function of $U$ [note that $F_U(1^-) = \mathbf{P}(X_1 < X_2, Y_1 > Y_2) = 1/4$]. We have

LEMMA 5.

$$\mathbf{P}(\Gamma) = \int_0^1 (1 - u)^{n-2} \, dF_U(u).$$

PROOF. Immediate since $U$ is determined by $P_1$ and $P_2$, and for all $k \geq 3$ the points $P_k$ are independent of $U$ and $\mathbf{P}(AX_k + BY_k > C | P_1, P_2) = U$ (the atom at $u = 1$ does not matter since the integrand vanishes there). $\square$

---

and $p'$ can never be comparable, i.e., neither $\mathrm{supp}(p) \subset \mathrm{supp}(p')$ nor $\mathrm{supp}(p) \supset \mathrm{supp}(p')$ can hold).

[9]$A, B$ and $C$ are thus random variables that are functions of $P_1$ and $P_2$.

Corollary 6.

$$\mathbf{P}(D_{12} = 1) = \mathbf{P}(\Gamma) = (n-2) \int_0^1 (1-u)^{n-3} F_U(u) \, du.$$

Proof. Integrate by parts:

$$\int_0^1 (1-u)^{n-2} \, dF_U(u) = [(1-u)^{n-2} F_U(u)]_0^1$$
$$+ (n-2) \int_0^1 (1-u)^{n-3} F_U(u) \, du,$$

and note that the first term vanishes. □

2.4. *Second moment.* To evaluate $\mathbf{P}(D_{12} = D_{13} = 1)$, we need the entries in the third column of the matrix $R$ as well. Let $Z_i = R(i, 3)$ be $n$ random variables that are $F$-distributed, with all the $X_i, Y_i, Z_i$ independent. Let $\Gamma'$ be the event that $D_{13} = 1$ (we will use $'$ for the $XZ$-problem), that is,

$$\Gamma' := \{X_1 < X_3, Z_1 > Z_3, A'X_k + B'Z_k < C' \text{ for all } k \neq 1, 3\},$$

where $A', B'$ and $C'$ are determined by $P_1' = (X_1, Z_1)$ and $P_3' = (X_3, Z_3)$ [cf. (1)]. Let $U'$ be the corresponding random variable: $U' := \mathbf{P}(A'X_2 + B'Y_2 > C'|P_1', P_3')$ if $X_1 < X_3$ and $Z_1 > Z_3$, and $U' := 1$ otherwise; put $W := \max\{U, U'\}$, with cumulative distribution function $F_W$.

Proposition 7.

$$\mathbf{P}(D_{12} = D_{13} = 1) = \mathbf{P}(\Gamma \cap \Gamma') \le (n-3) \int_0^1 (1-u)^{n-4} F_W(u) \, du.$$

Proof. For each $k \ge 4$ we have

$$\mathbf{P}(AX_i + BY_i < C, A'X_i + B'Z_i' < C'|P_1, P_2, P_1', P_3')$$
$$\le \min\{\mathbf{P}(AX_k + BY_k < C|P_1, P_2), \mathbf{P}(A'X_k + B'Z_k' < C'|P_1', P_3')\}$$
$$= \min\{1 - U, 1 - U'\} = 1 - \max\{U, U'\} = 1 - W.$$

Therefore

$$\mathbf{P}(\Gamma \cap \Gamma') \le \int_0^1 (1-u)^{n-3} \, dF_W(u).$$

As in Corollary 6, integrating by parts yields the result. □

2.5. *Poisson approximation.* The "Chen–Stein method" yields Poisson approximations for sums of Bernoulli random variables whose dependence is not too large. We will use the following formulation due to Arratia, Goldstein and Gordon [1]:

THEOREM 8. *Let $I$ be an arbitrary index set. For each $\alpha \in I$, let $Z_\alpha$ be a Bernoulli random variable with $\mathbf{P}(Z_\alpha = 1) = 1 - \mathbf{P}(Z_\alpha = 0) = p_\alpha > 0$, and let $B_\alpha \subset I$ be the "neighborhood of dependence" for $\alpha$; that is, $\alpha \in B_\alpha$ and $Z_\alpha$ is independent of $Z_\beta$ for all $\beta \notin B_\alpha$. Put*

$$Z := \sum_{\alpha \in I} Z_\alpha,$$

$$\lambda := \sum_{\alpha \in I} \mathbf{E}(Z_\alpha) = \sum_{\alpha \in I} p_\alpha,$$

$$b_1 := \sum_{\alpha \in I} \sum_{\beta \in B_\alpha} \mathbf{E}(Z_\alpha)\mathbf{E}(Z_\beta) = \sum_{\alpha \in I} \sum_{\beta \in B_\alpha} p_\alpha p_\beta,$$

$$b_2 := \sum_{\alpha \in I} \sum_{\beta \in B_\alpha \setminus \{\alpha\}} \mathbf{E}(Z_\alpha Z_\beta).$$

*Then*

$$\|\mathcal{L}(Z) - \mathrm{Poisson}(\lambda)\| \leq 2(b_1 + b_2)\frac{1 - e^{-\lambda}}{\lambda} \leq 2(b_1 + b_2).$$

PROOF. Theorem 1 in Arratia, Goldstein and Gordon [1], with no "near-independence" (i.e., $b_3' = b_3 = 0$).   □

2.6. *Notation.* We use the following standard notation, all as $n \to \infty$: $g(n) \sim h(n)$ for $\lim_n g(n)/h(n) = 1$; $g(n) \lesssim h(n)$ for $\limsup_n g(n)/h(n) \leq 1$; $g(n) \approx h(n)$ for $0 < \lim_n g(n)/h(n) < \infty$; $g(n) = O(h(n))$ for $\limsup_n g(n)/h(n) < \infty$; and $g(n) = o(h(n))$ for $\lim_n g(n)/h(n) = 0$. Also, log is the natural logarithm $\log_e$ throughout.

**3. The convex hull of $n$ random points in the plane.** Interestingly, the expectation $\mu_n$ of $S_2^{(n)}$ is related to the number of vertices, or edges, of the convex hull $K$ of the $n$ random points in the plane $P_1, P_2, \ldots, P_n$ (the connection does not, however, extend beyond the first moments). Denote that number by $V \equiv V^{(n)}$, and let $V_0$ be the number of edges of $K$ whose outward normal is positive.[10] The distribution $F$ is called *symmetric* if $F(-x) = 1 - F(x)$ for all $x$ [or, more generally, if there exists $x_0$ such that $F(x_0 - x) = 1 - F(x_0 + x)$ for all $x$].

---

[10]The "outward normal" to an edge of $K$ is perpendicular to the edge and points away from the interior of $K$.

PROPOSITION 9.

$$2\mu_n = \mathbf{E}(V_0) \geq \mathbf{P}(V_0 > 0) = 1 - \frac{1}{n}.$$

*Moreover, if $F$ is symmetric, then*

$$8\mu_n = \mathbf{E}(V).$$

PROOF. Let $E_{ij}$ be the indicator that the line segment $P_i P_j$ is an edge of $K$ with positive outward normal; then $V_0 = \sum_{i<j} E_{ij}$. Clearly, $\mathbf{P}(E_{ij} = 1) = \mathbf{P}(E_{12} = 1) = 2\mathbf{P}(\Gamma)$ (if the additional condition $X_1 < X_2, Y_1 > Y_2$ in $\Gamma$ is not satisfied, interchange $P_1$ and $P_2$; this yields the factor 2), and so $\mathbf{E}(V_0) = (n(n-1)/2)2\mathbf{P}(\Gamma) = 2\mu_n$.

Now $V_0 = 0$ if and only if there is a point $P_i$ that is maximal in both the $X$- and the $Y$-direction, that is, $X_i = \max_j X_j$ and also $Y_i = \max_j Y_j$. The probability of this event is $1/n$ (letting $i$ be the index where $X_i = \max_j X_j$, the probability that $Y_i = \max_j Y_j$ is $1/n$, since the $Y$'s are independent of the $X$'s). Therefore,

$$\mathbf{E}(V_0) \geq \mathbf{P}(V_0 \geq 1) = 1 - \frac{1}{n}.$$

If $F$ is symmetric, the same holds for outward normals in each of the four quadrants, and so $\mathbf{E}(V) = 4\mathbf{E}(V_0)$. □

Our main result for the number of vertices $V^{(n)}$ is:

THEOREM 10. *Let $F$ be a symmetric distribution. Then, as $n \to \infty$:*

   (i) *if $F \in \mathcal{EF}$, then $\mathbf{E}(V^{(n)}) \to \infty$; and*
   (ii) *if $F \in \mathcal{SE}$, then $\mathbf{E}(V^{(n)}) \to 4$ and $\mathbf{P}(V^{(n)} = 4) \to 1$.*

PROOF. Combine Proposition 9 above with results that will be obtained in the next two sections: Proposition 12 for (i) and Corollary 20 for (ii). □

Some intuition for the interesting result (ii) for "heavy-tailed" distributions is provided immediately after the proof of Theorem 19 in Section 5.2.[11] Figures[12] 1 and 2 show, for each one of five different distributions, $n = 10{,}000$ random points together with their convex hull and the resulting number of vertices $V^{(n)}$. In the context of random points drawn from radially symmetric distributions (rather than independent coordinates), Carnal [4] has shown that $\mathbf{E}(V^{(n)})$ converges to a constant $\geq 4$ for a certain class of heavy-tailed distributions (with the constant depending on the distribution).

---

[11] Fisher [8] shows that for certain distributions (including the Weibull distributions with parameter $0 < \alpha < 1$) the limit shape of the normalized convex hull is $\{(x, y) \in \mathbb{R}^2 : |x| + |y| \leq 1\}$—which is the convex hull of four points. However, this does *not* imply that the number of vertices $V^{(n)}$ converges to 4, since there may be many vertices close to each one of these four points (as is the case for the uniform distribution, where the limit shape is the unit square, and $V^{(n)} \to \infty$).

[12] Generated by MAPLE.

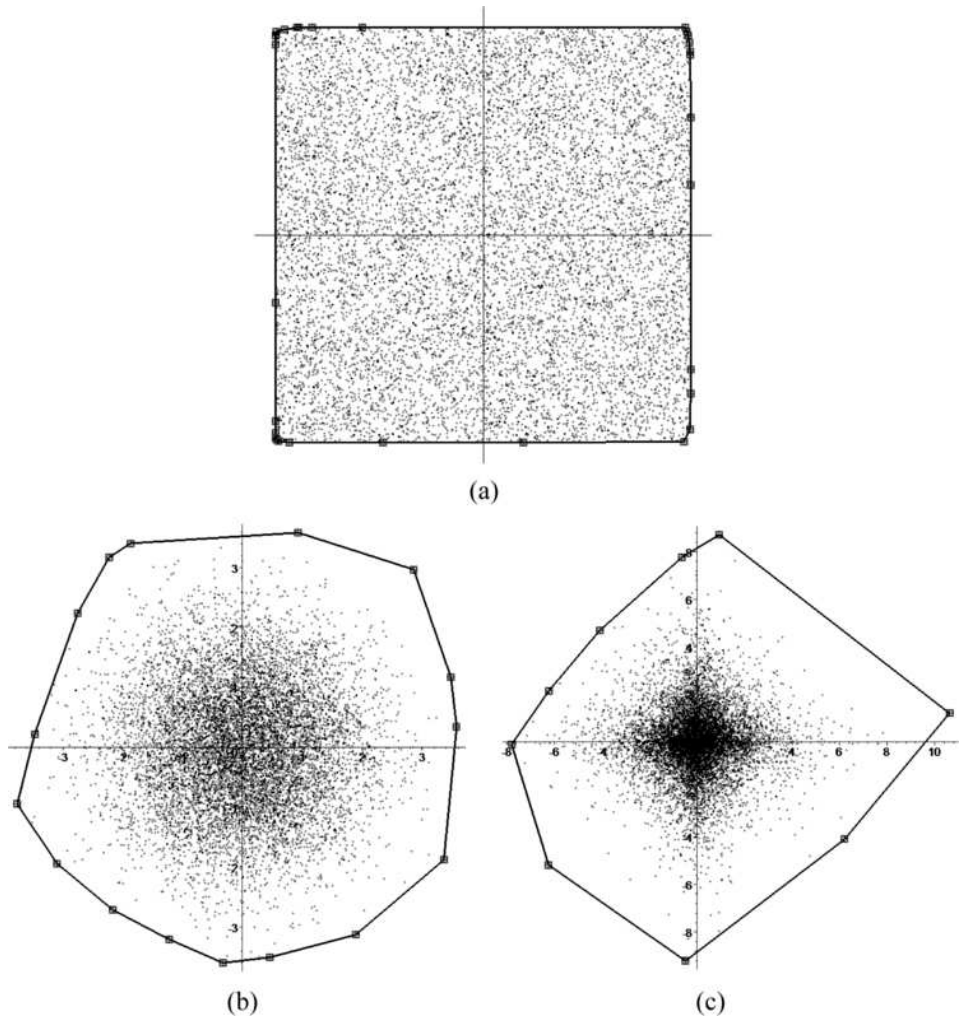

FIG. 1.   *The number of vertices V of the convex hull of n random points drawn from three distributions in $\mathcal{EF}$.* (a) *Uniform distribution:* $n = 10{,}000$, $V = 29$. (b) *Normal distribution:* $n = 10{,}000$, $V = 16$. (c) *Exponential distribution:* $n = 10{,}000$, $V = 9$.

We conclude this section with a lemma that is useful when comparing distributions (see its use in the next section).

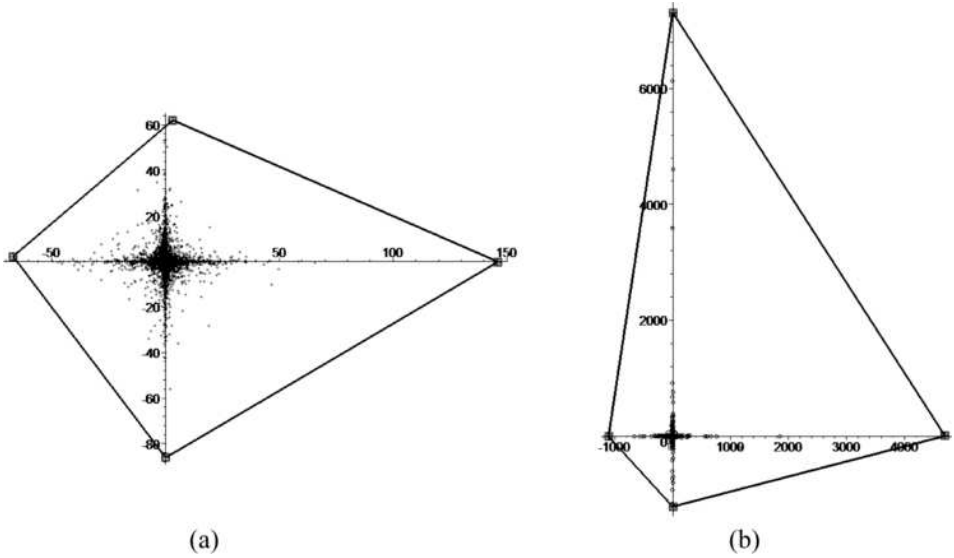

**(a)** **(b)**

Fig. 2. *The number of vertices $V$ of the convex hull of $n$ random points drawn from two distributions in $\mathcal{SE}$.* (a) *Weibull distribution, $\alpha = 1/2$: $n = 10{,}000$, $V = 4$.* (b) *Cauchy distribution: $n = 10{,}000$, $V = 4$.*

LEMMA 11. *Let $F'$ and $F''$ be two distributions, with supports[13] $(a', b')$ and $(a'', b'')$ and corresponding $\mu'_n$ and $\mu''_n$. If there exists a strictly increasing convex function $\varphi : (a'', b'') \to (a', b')$ such that $F''(x) = F'(\varphi(x))$ for all $x \in (a'', b'')$, then $\mu'_n \leq \mu''_n$.*

PROOF. Let $(X'_i)_{1 \leq i \leq n}$ and $(Y'_i)_{1 \leq i \leq n}$ be independent and $F'$-distributed random variables, and define $X''_i := \varphi^{-1}(X'_i)$ and $Y''_i := \varphi^{-1}(Y'_i)$. Put $P'_i = (X'_i, Y'_i)$ and $P''_i = (X''_i, Y''_i)$, and let $K'$ and $K''$ be the convex hulls of $\{P'_i\}_i$ and $\{P''_i\}_i$, respectively. Since $\mathbf{P}(X''_i \leq x) = \mathbf{P}(\varphi^{-1}(X'_i) \leq x) = \mathbf{P}(X'_i \leq \varphi(x)) = F'(\varphi(x)) = F''(x)$, the $(X''_i)_i$ and $(Y''_i)_i$ are $F''$-distributed. If $Ax + By$ is a supporting line to $K'$ at $P'_i$ with $A, B > 0$, then $(Ap)x + (Bq)y$ is a supporting line to $K''$ at $P''_i$, where $p, q > 0$ are subgradients of $\varphi$ at $X'_i$ and $Y'_i$, respectively. Therefore $V'_0 + 1 \leq V''_0 + 1$ (the number of vertices supported by positive outward normals is larger by one than the number of edges supported by such normals), and so $\mu'_n \leq \mu''_n$. $\quad\square$

## 4. Exponential and faster decreasing tails.

### 4.1. *The class $\mathcal{EF}$.* We define the class of distributions $\mathcal{EF}$ with "*Exponential and Faster decreasing tails*" as those continuous distributions $F$

---

[13] $-\infty \leq a' < b' \leq \infty$ and $-\infty \leq a'' < b'' \leq \infty$.

with support $(a, b)$ (where $-\infty \le a < b \le \infty$) whose "tail" $G(x) = 1 - F(x)$ is a log-concave function; that is, $G(x) = e^{-g(x)}$ where $g : (a, b) \to (0, \infty)$ is a strictly increasing convex function. The functions $G$ and $g = -\log G$ are usually called the *survival function* and the *cumulative hazard function*, respectively; for a collection of results on log-concave probabilities, see Bagnoli and Bergstrom [2].[14] A sufficient (but not necessary) condition for the log-concavity of $G$ is that the density function $f = F'$ be continuously differentiable and log-concave. Some distributions included in the class $\mathcal{EF}$ are the following (for simplicity, we take standard normalizations; replacing $x$ with $\lambda x + \nu$ for any $\lambda > 0$ and $\nu$ clearly preserves the log-concavity of $G$):

- *Exponential*: $G(x) = e^{-x}$ for $x \in (0, \infty)$.
- *Normal*: $G(x) = \int_x^\infty (2\pi)^{-1/2} e^{-y^2/2} \, dy$ for $x \in (-\infty, \infty)$.
- *Weibull with parameter* $\alpha \ge 1$: $G(x) = e^{-x^\alpha}$ for $x \in (0, \infty)$, where $\alpha \ge 1$ (these are the $F_\alpha$ of the Introduction).
- $G(x) = e^{-e^x}$ for $x \in (-\infty, \infty)$.
- *Logistic*: $G(x) = 1/(1 + e^x)$ for $x \in (-\infty, \infty)$.
- *Uniform*: $G(x) = 1 - x$ for $x \in (0, 1)$.

Each such distribution is by definition an increasing convex transformation of the exponential distribution: if $F(x) = 1 - e^{-g(x)}$, then $F(x) = F^{\exp}(g(x))$ for every $x$ in the support of $F$ [where $F^{\exp}(x) = 1 - e^{-x}$ is the exponential cumulative distribution function]. By Lemma 11, it thus follows that the exponential distribution yields the lower bound on $\mu_n$ over the class $\mathcal{EF}$. Now Haigh [11] proved that $\mu_n^{\exp} \approx \log\log n$, and so we have

PROPOSITION 12. *If $F \in \mathcal{EF}$, then $\mu_n \to \infty$ as $n \to \infty$.*

PROOF. If $F \in \mathcal{EF}$, then $\mu_n^F \ge \mu_n^{\exp} \approx \log\log n \to \infty$ by Lemma 11 and Haigh [11]. □

Rényi and Sulanke [13] provide more precise results: $\mu_n^{\text{normal}} \approx \sqrt{\log n}$ and $\mu_n^{\text{uniform}} \approx \log n$. Also, we note that the class $\mathcal{EF}$ can be taken to be much larger; see Section 6(b).

4.2. *Poisson approximation.* Our Theorem 1 for the class $\mathcal{EF}$ is an immediate consequence of Proposition 12, together with the general result of Theorem 13 below (which holds for *any* distribution $F$, not necessarily in $\mathcal{EF}$). The analysis will also yield the universal upper bound of Theorem 17.

---

[14]The class of positive random variables with a log-concave $G$ is usually called *IFR* (for *Increasing Failure Rate*).

THEOREM 13. *For every distribution $F$,*

$$\|\mathcal{L}(S_2^{(n)}) - \text{Poisson}(\mu_n)\| = O\left(\frac{1}{\sqrt{\mu_n}}\right) \qquad \text{as } n \to \infty.$$

The remainder of this section is devoted to the proof of Theorem 13. For every $x \in \mathbb{R}$ and $u \in (0, 1)$, let $\nu(x; u) := \mathbf{P}(U < u | X_1 = x)$ [recall the definition (3) of $U$]; then

(4)
$$F_U(u) = \int_{-\infty}^{\infty} \nu(x; u) \, dF(x)$$

and

(5)
$$F_W(u) = \int_{-\infty}^{\infty} \nu(x; u)^2 \, dF(x),$$

since, given $X_1 = x$, the events $U < u$ and $U' < u$ are independent (the first depends on $Y_1, X_2, Y_2$ and the second on $Z_1, X_3, Z_3$; see Section 2.4).

For every $b \geq 0$ and $u \in (0, 1)$, let $\kappa_u(b)$ be determined by the equation $\mathbf{P}(X + bY \geq \kappa_u(b)) = u$ (it is unique since $X + bY$ is a continuous random variable). Let $K \equiv K_u := \{(x, y) \in \mathbb{R}^2 : x + by < \kappa_u(b) \text{ for all } b > 0\}$ be the set of all points that are not contained in any half-plane of probability $u$ with positive normal (see Figure 3). Clearly, if either $P_1 \in K$ or $P_2 \in K$, then $U \geq u$ [since, for all $b > 0$, the line $x + by = c$ through that point has $c < \kappa_u(b)$ and so $\mathbf{P}(X + bY > c) \geq \mathbf{P}(X + bY \geq \kappa_u(b)) = u$]. The set $K$ is convex (it is an intersection of half-spaces) and comprehensive [i.e., $(x', y') \leq (x, y) \in K$ implies that $(x', y') \in K$]. Let $y = \eta(x; u)$ be the equation of its boundary, that is, $\eta(x; u) := \sup\{y : (x, y) \in K\}$ [with $\eta(x; u) := -\infty$ when there is no such $y$]. We have:

LEMMA 14. *For every $x$ and $u \in (0, 1)$*

$$\nu(x; u) \leq u G(\eta(x; u)) \leq u.$$

PROOF. Let $P_1 = (x_1, y_1)$. If $P_1 \in K$, then, as we saw above, $\mathbf{P}(U < u | X_1 = x_1, Y_1 = y_1) = 0$.

If $P_1 \notin K$, then $y_1 \geq \eta(x_1; u)$ (again, see Figure 3); let $b_0 \equiv b_0(x_1, y_1) := \inf\{b > 0 : x_1 + by_1 \geq \kappa_u(b)\}$. The function $\kappa_u$ is continuous since the distribution $F$ is continuous, and so $x_1 + b_0 y_1 \geq \kappa_u(b_0)$ [note that we may well have $b_0 = 0$, for which $\kappa_u(0) = G^{-1}(u)$]. Assume that $U < u$; then there exists $b > 0$ such that $X_2 + bY_2 = x_1 + by_1 \geq \kappa_u(b)$, and so $b \geq b_0$. Now $Y_2 < y_1$ (since $U < u \leq 1$); therefore $X_2 + b_0 Y_2 \geq x_1 + b_0 y_1$, which, as we saw above, is $\geq \kappa_u(b_0)$. Thus $U < u$ implies that $P_2$ lies above the line $x + b_0 y = \kappa_u(b_0)$, and so $\mathbf{P}(U < u | X_1 = x_1, Y_1 = y_1) \leq \mathbf{P}(X_2 + b_0 Y_2 \geq \kappa_u(b_0)) \leq u$ by definition of $\kappa_u$.

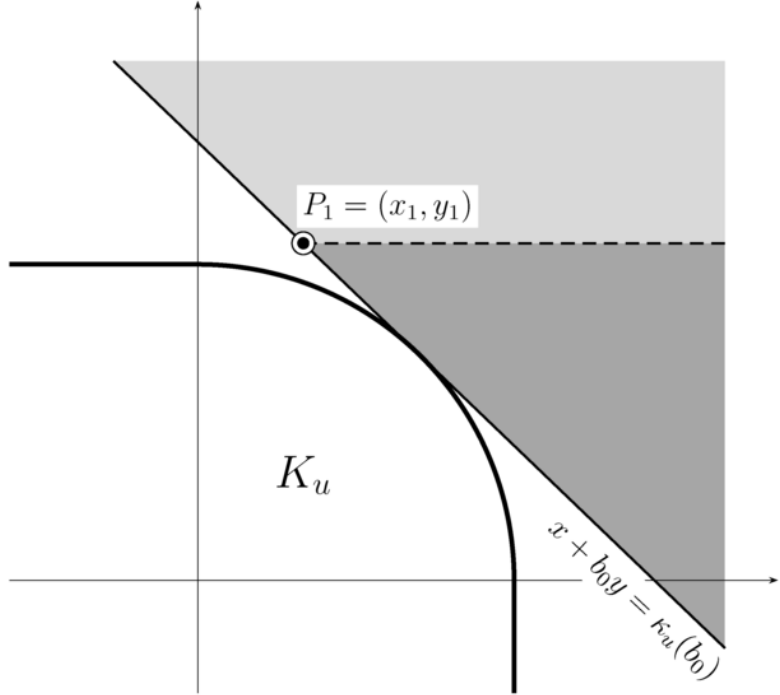

Fig. 3. *If $U < u$, then $P_2$ lies in the darkly shaded area; the probability of the whole shaded area is $u$ (Lemma 14).*

Taking expectation over $Y_1 = y_1$ therefore yields

$$\mathbf{P}(U < u | X_1 = x_1) \leq 0\mathbf{P}(P_1 \in K | X_1 = x_1) + u\mathbf{P}(P_1 \notin K | X_1 = x_1)$$

$$\leq u\mathbf{P}(Y_1 \geq \eta(x_1; u)) = uG(\eta(x_1; u)) \leq u. \qquad \square$$

LEMMA 15. *For every $x$ and $u \in (0, 1)$*

$$G(x)G(\eta(x; u)) \leq u.$$

PROOF. If $P_1 = (x_1, y_1) \notin K$, then (see Figure 4) there exists $b > 0$ such that $c := x_1 + by_1 \geq \kappa_u(b)$, and so $\mathbf{P}(X + bY \geq c) \leq u$. Therefore, $\mathbf{P}(X \geq x_1, Y \geq y_1) \leq \mathbf{P}(X + bY \geq c) \leq u$, and so $G(x_1)G(y_1) \leq u$. This holds for all $y_1 > \eta(x_1; u)$, and $G$ is a continuous function. $\square$

Combining the inequalities in the last two lemmas yields:

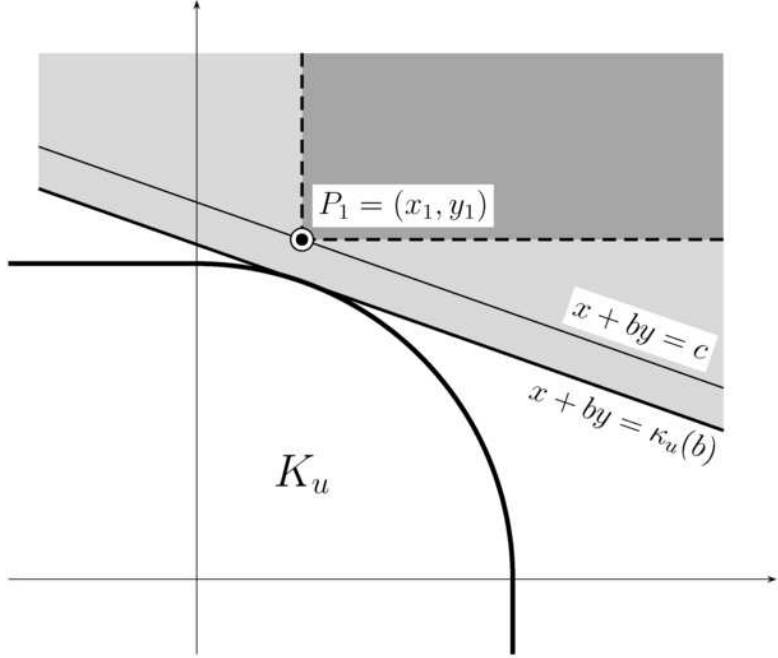

FIG. 4. *The probability of the whole shaded area is* $u$, *and that of the darkly shaded area is* $G(x_1)G(y_1)$ *(Lemma 15).*

COROLLARY 16. *For every* $x$ *and* $u \in (0,1)$

$$\nu(x;u) \leq \min\left\{u, \frac{u^2}{G(x)}\right\}.$$

From this we can immediately obtain an upper bound on $\mu_n$ which applies to any distribution $F$. This bound is known; see Devroye [6].

COROLLARY 17. *For every distribution* $F$,

$$\mathbf{P}(\Gamma) = O\left(\frac{\log n}{n^2}\right) \quad and \quad \mu_n = O(\log n) \qquad as \ n \to \infty.$$

PROOF. Let $t$ be such that $G(t) = u$. Applying Corollary 16 in the formula (4) yields

$$F_U(u) = \int_{-\infty}^{\infty} \nu(x;u)\,dF(x) \leq \int_{-\infty}^{t} \frac{u^2}{G(x)}\,dF(x) + \int_{t}^{\infty} u\,dF(x)$$

$$= u^2 \int_u^1 \frac{1}{z} \, dz + u^2 = u^2 \left( \log \frac{1}{u} + 1 \right)$$

[we have used the substitution $z = G(x)$]. Therefore,

$$\mathbf{P}(\Gamma) = (n-2) \int_0^1 (1-u)^{n-3} F_U(u) \, du$$

$$\leq (n-2) \int_0^1 (1-u)^{n-3} u^2 \left( \log \frac{1}{u} + 1 \right) du$$

$$\leq \frac{2 \log n}{n^2} + O\left( \frac{1}{n^2} \right). \qquad \square$$

We can now prove Theorem 13.

PROOF OF THEOREM 13.   Let $\xi(u)$ be determined by $G(\xi(u)) = \sqrt{F_U(u)}$. For $x \leq \xi(u)$, we will use the inequality $\nu(x; u) \leq u^2/G(x)$ to get

$$\int_{-\infty}^{\xi(u)} \nu(x; u)^2 \, dF(x) \leq \int_{-\infty}^{\xi(u)} \nu(x; u) \frac{u^2}{G(x)} \, dF(x)$$

$$\leq \frac{u^2}{\sqrt{F_U(u)}} \int_{-\infty}^{\xi(u)} \nu(x; u) \, dF(x)$$

$$\leq \frac{u^2}{\sqrt{F_U(u)}} \int_{-\infty}^{\infty} \nu(x; u) \, dF(x)$$

$$= \frac{u^2}{\sqrt{F_U(u)}} F_U(u) = u^2 \sqrt{F_U(u)}.$$

For $x \geq \xi(u)$, we use $\nu(x; u) \leq u$ to get

$$\int_{\xi(u)}^{\infty} \nu(x; u)^2 \, dF(x) \leq \int_{\xi(u)}^{\infty} u^2 \, dF(x) = u^2 G(\xi(u)) = u^2 \sqrt{F_U(u)}.$$

Altogether,

$$F_W(u) \leq 2u^2 \sqrt{F_U(u)}.$$

Therefore,

$$\int_0^1 (1-u)^{n-4} F_W(u) \, du$$

$$\leq 2 \int_0^1 (1-u)^{n-4} u^2 \sqrt{F_U(u)} \, du$$

$$\leq 2 \left( \int_0^1 (1-u)^{n-3} F_U(u) \, du \right)^{1/2} \left( \int_0^1 (1-u)^{n-5} u^4 \, du \right)^{1/2},$$

by the Cauchy–Schwarz inequality. The first integral is $\mathbf{P}(\Gamma)/(n-2) = O(\mu_n n^{-3})$ [by Corollary 6 and (2)], the second integral is $O(n^{-5})$, and so

$$\int_0^1 (1-u)^{n-4} F_W(u)\, du = O(\mu_n^{1/2} n^{-4}).$$

Therefore, by Proposition 7,

$$\mathbf{P}(\Gamma \cap \Gamma') = O(\mu_n^{1/2} n^{-3}).$$

We now apply Theorem 8 to $S_2^{(n)} = \sum_{i<j} D_{ij}$. There are $n(n-1)/2 = O(n^2)$ terms $D_{ij}$; the neighborhood of dependence of each $D_{ij}$ consists of $D_{ik}$ and $D_{jk}$ for all $k$, and so it is of size $2n - 3 = O(n)$. Therefore,

$$b_1 = O(n^2)O(n)\mathbf{P}(\Gamma)^2 = O(n^3(\mu_n n^{-2})^2) = O(\mu_n^2 n^{-1})$$

and

$$b_2 = O(n^2)O(n)\mathbf{P}(\Gamma \cap \Gamma') = O(\mu_n^{1/2}).$$

This yields $\|\mathcal{L}(S_2^{(n)}) - \text{Poisson}(\mu_n)\| \leq 2(b_1 + b_2)/\mu_n = O(\mu_n n^{-1} + \mu_n^{-1/2})$. Now $\mathbf{E}(V^{(n)}) = O(\log n)$, and so $\mu_n = O(\log n)$, for any distribution $F$; this follows from Theorem 1 (for dimension $d = 2$) in [6]. Therefore $\|\mathcal{L}(S_2^{(n)}) - \text{Poisson}(\mu_n)\| = O(\mu_n^{-1/2})$. $\quad\square$

PROOF OF THEOREM 1. Combine Proposition 12 and Theorem 13. $\quad\square$

## 5. Slower than exponential decreasing tails.

5.1. *The class $\mathcal{SE}$.* We define the class of distributions $\mathcal{SE}$ with "*Slower than Exponential decreasing tails*" as those distributions $F$ with support $(a, \infty)$ (where $a \geq -\infty$) whose tail $G = 1 - F$ satisfies the following two conditions:

[SE1] "*Subexponentiality*":

$$\mathbf{P}(X_+ + Y_+ > t) \sim 2G(t) \qquad \text{as } t \to \infty,$$

where $X, Y$ are independent $F$-distributed random variables and $Z_+ := \max\{Z, 0\}$; and

[SE2] "*Uniformity*": For all $c > 1$,

(6) $$G(ct) \gtrsim G(t)^c \qquad \text{as } t \to \infty, \qquad \text{uniformly as } c \to 1^+;$$

that is, for every $\varepsilon > 0$ there exist $t_0 \equiv t_0(\varepsilon)$ and $c_0 \equiv c_0(\varepsilon) > 1$ such that $G(ct)/G(t)^c > 1 - \varepsilon$ for all $t > t_0$ and all $c \in (1, c_0)$.

Distributions satisfying [SE1] are called *subexponential distributions* (originally introduced by Chistyakov [5]). Some examples are (see Table 3.7 in Goldie and Klüppelberg [9]; again, we use standard normalizations for simplicity):

- *Regularly varying tails*: $G(x) = x^{-\alpha}\ell(x)$, where $\alpha \geq 0$ and $\ell$ is a *slowly varying function*, that is, $\lim_{x\to\infty}\ell(cx)/\ell(x) = 1$ for every $c > 0$. This includes:
  - *Pareto*: $G(x) = x^{-\alpha}$ for $x \in (1, \infty)$, where $\alpha > 0$.
  - *Cauchy*: $G(x) = \int_x^\infty (\pi(1 + y^2))^{-1}\,dy = \arctan(x)/\pi + 1/2$ for $x \in (0, \infty)$.
  - *$\alpha$-stable*, where $0 < \alpha < 2$.
- *Lognormal*: $G(x) = \int_x^\infty (\sqrt{2\pi}y)^{-1}e^{-\log^2 y/2}\,dy$ for $x \in (0, \infty)$.
- *Weibull with parameter* $0 < \alpha < 1$: $G(x) = e^{-x^\alpha}$ for $x \in (0, \infty)$.
- "*Almost*" *exponential*: $G(x) = e^{-x(\ln x)^{-\alpha}}$ for $x \in (1, \infty)$, where $\alpha > 0$.

(However, the exponential distribution does *not* satisfy [SE1].)

As for condition [SE2], in terms of the cumulative hazard function $g(t) := -\log G(t)$, it says that for every $\varepsilon > 0$ there exist $t_0 \equiv t_0(\varepsilon)$ and $c_0 \equiv c_0(\varepsilon) > 1$ such that $g(ct) \leq cg(t) + \varepsilon$ for all $t > t_0$ and all $c \in (1, c_0)$. Therefore, a sufficient condition for [SE2] is that $g(t)/t$ be a nonincreasing function for large enough[15] $t$; this is the case when $g$ is concave (and so $G$ is *log-convex*; contrast with $\mathcal{EF}$), or even star-concave[16] (we will see in Lemma 18(ii) below that [SE1] implies $g(t)/t \to 0$ as $t \to \infty$). It is now easy to verify that all the distributions listed above also satisfy [SE2]. Finally, $\mathcal{SE}$ is closed under "tail equivalence": if $1 - F(t) \sim 1 - F'(t)$ as $t \to \infty$, then $F \in \mathcal{SE}$ if and only if $F' \in \mathcal{SE}$ (for [SE1], see Theorem 3 in Teugels [14]).

The next lemma collects a number of properties that will be used in the proof below.

Lemma 18. *Let $F$ satisfy* [SE1]. *Then:*

(i)

$$\mathbf{P}(X + Y > t) \lesssim 2G(t) \qquad as\ t \to \infty. \tag{7}$$

(ii) $g(t) := -\log G(t) = o(t)$ *as* $t \to \infty$.

(iii) *There exist* $\gamma_t > 0$ *such that*

$$\lim_{t\to\infty}\gamma_t = 0, \tag{8}$$

$$\lim_{t\to\infty}\gamma_t t = \infty, \tag{9}$$

$$\lim_{t\to\infty}G(t)^{\gamma_t} = 1. \tag{10}$$

---

[15] The class of positive random variables where $g(t)/t$ is a nonincreasing function for *all* $t$ is usually called DFRA (for *Decreasing Failure Rate Average*).

[16] That is, $g(\lambda x) \geq \lambda g(x) + (1 - \lambda)g(0)$ for all $x \geq 0$ and all $0 \leq \lambda \leq 1$.

*Moreover, if F also satisfies* [SE2], *then*

$$\lim_{t \to \infty} \frac{G((1 + \gamma_t)t)}{G(t)} = 1. \tag{11}$$

PROOF. (i) is immediate from [SE1] since $X + Y \leq X_+ + Y_+$. As for (ii), it is a well-known property of subexponential distributions (e.g., it follows from (1.4) in Goldie and Klüppelberg [9]). To get (iii), take, for example, $\gamma_t = 1/\sqrt{tg(t)}$, and then (8), (9) and (10) immediately follow from (ii); finally, (10) together with (6) imply (11). □

5.2. *First moment.* In this section we will prove that, for distributions in $\mathcal{SE}$, the expected number of two-point ESS converges to 1/2, and the number of vertices of the convex hull converges in probability to 4. Some intuition is provided after the proof of Theorem 19. The main result is

THEOREM 19. *Let* $F \in \mathcal{SE}$. *Then*

$$\mathbf{P}(\Gamma) \sim \frac{1}{n^2} \quad and \quad \mu_n \to \frac{1}{2} \qquad as \ n \to \infty.$$

As a result, the number of vertices $V^{(n)}$ of the convex hull of $n$ random points satisfies

COROLLARY 20. *Let* $F$ *be a symmetric* $\mathcal{SE}$ *distribution. Then*

$$\mathbf{E}(V^{(n)}) \to 4 \quad and \quad \mathbf{P}(V^{(n)} = 4) \to 1 \qquad as \ n \to \infty.$$

PROOF. Theorem 19 and Proposition 9 yield $\mathbf{E}(V_0) \to 1$ and $\mathbf{P}(V_0 = 0) = 1/n \to 0$, and so $\mathbf{P}(V_0 \neq 1) \to 0$. The result follows since $V = 4V_0$. □

Thus, for symmetric $\mathcal{SE}$ distributions, the probability that the convex hull is a quadrilateral converges to 1.

*For the remainder of this section we assume that* $F \in \mathcal{SE}$.

The proof of Theorem 19 uses the following result:

PROPOSITION 21. *As* $u \to 0$

$$F_U(u) \lesssim \tfrac{1}{2} u^2.$$

Before proving Proposition 21 (to which most of this section is devoted), we use it to prove Theorem 19.

Proof of Theorem 19.   Given $\varepsilon > 0$, let $\delta > 0$ be such that $F_U(u) \leq (1 + \varepsilon)u^2/2$ for all $u < \delta$. We will use Corollary 6, and separate the integral into two parts. For the first part, we have

$$\int_0^\delta (1-u)^{n-3} F_U(u)\, du \leq (1+\varepsilon)\frac{1}{2} \int_0^\delta (1-u)^{n-3} u^2 \, du$$

$$\leq (1+\varepsilon)\frac{1}{2} \int_0^1 (1-u)^{n-3} u^2 \, du$$

$$= (1+\varepsilon)\frac{1}{2}\frac{2}{n(n-1)(n-2)}.$$

As for the second part, we get

$$\int_\delta^1 (1-u)^{n-3} F_U(u)\, du \leq (1-\delta)^{n-3},$$

which is less than, say, $\varepsilon/n^3$ for all $n$ large enough. Adding the two bounds, multiplying by $n-2$, and recalling Corollary 6 yields $\mathbf{P}(\Gamma) \leq (1+2\varepsilon)/n^2$ for all $n$ large enough. The opposite inequality is in Proposition 9 [recall (2)]. $\square$

The proof of Proposition 21 requires careful analysis. To get some intuition, consider the convex hull of $n$ random points $P_1, P_2, \ldots, P_n$. Let $P_i = (X_i, Y_i)$ be the (a.s. unique) point with maximal $X$-coordinate, that is, $X_i = \max_k X_k$. An essential property of subexponential distributions is that $X_i$ is much larger than all the other $X_k$ for $k \neq i$. In addition, the corresponding $Y$-coordinate, namely $Y_i$, is also much smaller than $X_i$. The same holds for the point $P_j = (X_j, Y_j)$ with maximal $Y$-coordinate, which implies that, with high probability, all the points $P_k$ with $k \neq i, j$ will lie well below the line connecting $P_i$ and $P_j$, so that $P_i$ and $P_j$ will be the only vertices with positive outward normals. This basic picture can be seen in Figure 5 (recall also Figure 2). The points in the region $L^2$ have large $X$ (bigger than an appropriate $t$), whereas the width of $L^2$ (in the $Y$-coordinate) is small relative to $t$. The same holds for the region $L^1$, with $X$ and $Y$ interchanged. These two regions will thus "catch," with high probability, the points $P_i$ and $P_j$ with maximal $X$ and maximal $Y$, respectively.

Fix $0 < \varepsilon < 1$, and let $t \equiv t_{u,\varepsilon}$ be such that

$$(12) \qquad\qquad\qquad G(t) = (1+\varepsilon)u;$$

then $u \to 0$ is equivalent to $t \to \infty$ (since $\varepsilon > 0$ is fixed). We will say that $t$ and $u$ correspond to one another if they are related by (12). Next, we define the following sets in $\mathbb{R}^2$ (see Figure 5):

$$L^0 \equiv L^0_{u,\varepsilon} := \{(x,y) : x \leq t, y \leq t, x + y \leq t\},$$

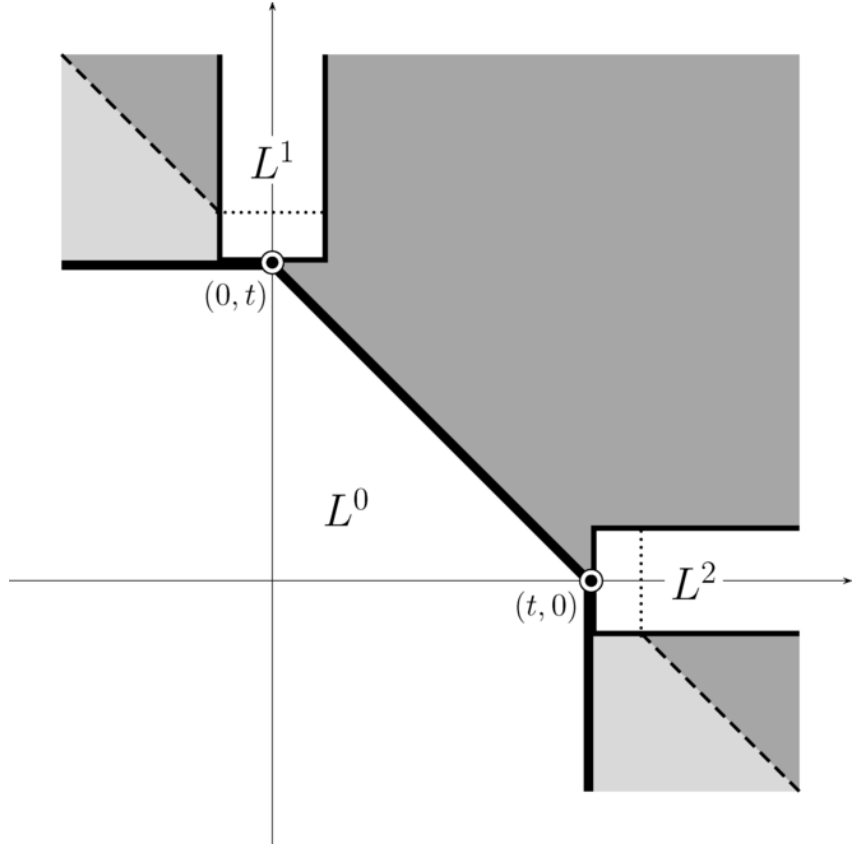

Fig. 5. $L^3$ is the whole shaded area and $L^4$ (see Lemma 29) is the darkly shaded area.

$$L \equiv L_{u,\varepsilon} := \mathbb{R}^2 \backslash L^0,$$
$$L^1 \equiv L^1_{u,\varepsilon} := \{(x,y) : |x| \leq \gamma_t t, y > t\},$$
$$L^2 \equiv L^2_{u,\varepsilon} := \{(x,y) : x > t, |y| \leq \gamma_t t\},$$
$$L^3 \equiv L^3_{u,\varepsilon} := L \backslash (L^1 \cup L^2).$$

The reader should keep in mind that $t$ as well as all the sets $L, L^0, \dots$ depend on $u$ (and $\varepsilon$).

For simplicity, we will write $\mathbf{P}(L)$ instead of $\mathbf{P}(P_i \in L)$.

LEMMA 22. *As $u \to 0$,*

$$\mathbf{P}(L_{u,\varepsilon}) = O(u).$$

PROOF.

$$\mathbf{P}(L) \leq \mathbf{P}(x+y>t) + \mathbf{P}(x>t) + \mathbf{P}(y>t) \lesssim 4G(t) = 4(1+\varepsilon)u.$$

□

LEMMA 23.   *There exists $u_0 \equiv u_0(\varepsilon) \in (0,1)$ such that, for all $u < u_0$, if $a, b, c$ satisfy $a, b > 0$ and $\mathbf{P}(aX + bY > c) < u$ (where $X, Y$ are independent and $F$-distributed), then $c/a > t_{u,\varepsilon}$ and $c/b > t_{u,\varepsilon}$.*

PROOF.   Assume without loss of generality that $a \geq b$. If $c/a \leq t \equiv t_{u,\varepsilon}$, then

$$\mathbf{P}(aX + bY > c) \geq \mathbf{P}(aX + bY > at)$$
$$\geq \mathbf{P}(|Y| \leq \gamma_t t, X > (1 + \gamma_t)t)$$
$$= H(\gamma_t t)G((1 + \gamma_t)t),$$

where $H(z) := 1 - F(-z) - G(z)$ for $z \geq 0$ (we used $a \geq b$ in the second inequality). Now $H(\gamma_t t)G((1 + \gamma_t)t)/G(t) \to 1$ as $t \to \infty$ by (9) and (11), and so $H(\gamma_t t)G((1 + \gamma_t)t) > G(t)/(1 + \varepsilon) = u$ for all $t$ large enough, or all $u$ small enough. This contradiction shows that indeed $c > at \geq bt$.   □

COROLLARY 24.   *For all $u < u_0$,*
$$\mathbf{P}(P_1 \in L^0_{u,\varepsilon} \text{ or } P_2 \in L^0_{u,\varepsilon}, U < u) = 0.$$

PROOF.   If $U < u < u_0$, then the entire set $L^0$ lies *below* the line $Ax + By = C$ through $P_1$ and $P_2$ [this holds for its two extreme points, $(t, 0)$ and $(0, t)$, by Lemma 23, and $A, B > 0$]; therefore $U < u$ implies that $P_1 \notin L^0$ and $P_2 \notin L^0$.   □

At this point we immediately get the following bounds:

PROPOSITION 25.
$$\mathbf{P}(U < u) = O(u^2) \qquad \text{as } u \to 0,$$
$$\mu_n = O(1) \qquad \text{as } n \to \infty.$$

PROOF.   Corollary 24 and Lemma 22 imply that $\mathbf{P}(U < u) \leq \mathbf{P}(P_1, P_2 \in L) = \mathbf{P}(L)^2 = O(u^2)$ as $u \to 0$. Using this in the computation of the proof of Theorem 19 yields $\mathbf{P}(\Gamma) = O(1/n^2)$, and so $\mu_n = O(1)$.   □

To get $\mu_n \lesssim 1/2$ will require a more refined analysis (the best constant we can get up to this point is $\mu_n \lesssim 4$). We start with a useful inequality:

LEMMA 26.   *Let $X$ and $Y$ be independent and $F$-distributed. Then, for every $a, b, c, \theta > 0$,*

$$(13) \quad \mathbf{P}(aX + bY > c) \geq H\left(\theta\frac{c}{b}\right)G\left((1 + \theta)\frac{c}{a}\right) + H\left(\theta\frac{c}{a}\right)G\left((1 + \theta)\frac{c}{b}\right).$$

Proof. We have

$$\mathbf{P}(aX + bY > c) \geq \mathbf{P}(aX > (1+\theta)c, |bY| \leq \theta c)$$
$$+ \mathbf{P}(|aX| \leq \theta c, bY > (1+\theta)c). \qquad \square$$

The next three lemmas will deal, respectively, with the three cases: (i) $P_1 \in L^1$ and $P_2 \in L^2$; (ii) $P_1, P_2 \in L^1$ or $P_1, P_2 \in L^2$; and (iii) $P_1, P_2 \in L^3$ (recall Corollary 24). The corresponding probabilities turn out to be of the order of $u^2/2$ in the first case, and $o(u^2)$ in the other two cases.

Lemma 27. As[17] $u \to 0$,

$$\mathbf{P}(P_1 \in L^1_{u,\varepsilon}, P_2 \in L^2_{u,\varepsilon}, U < u) \leq \tfrac{1}{2}u^2 + \varepsilon O(u^2).$$

Proof. Let $P_1 \in L^1$ and $P_2 \in L^2$ be such that $U < u < u_0$, where $u_0$ is given by Lemma 23, and let $t_0$ correspond to $u_0$. The line $P_1 P_2$ is $Ax + By = C$ with $A = Y_1 - Y_2$, $B = X_2 - X_1$ and $C = X_2 Y_1 - X_1 Y_2$. Since $U < u < u_0$, we have $C/A > t$ and $C/B > t$ by Lemma 23 and so, taking $\theta = \gamma_t$ in Lemma 26,

$$U = \mathbf{P}(AX + BY > C | P_1, P_2) \geq H(\gamma_t t)\Big[ G\Big((1+\gamma_t)\frac{C}{A}\Big) + G\Big((1+\gamma_t)\frac{C}{B}\Big) \Big].$$

Now $X_2, Y_1 > t$ and $|X_1|, |Y_2| < \gamma_t t$, and so $|X_1| < \gamma_t X_2$ and $|Y_2| < \gamma_t Y_1$, which implies that $A \geq Y_1(1 - \gamma_t), B \geq X_2(1 - \gamma_t)$ and $C \leq X_2 Y_1(1 + \gamma_t^2)$. Therefore,

$$\frac{C}{A} \leq X_2 \frac{1 + \gamma_t^2}{1 - \gamma_t} \quad \text{and} \quad \frac{C}{B} \leq Y_1 \frac{1 + \gamma_t^2}{1 - \gamma_t}.$$

Put $\rho_t := (1 + \gamma_t)(1 + \gamma_t^2)/(1 - \gamma_t) > 1$; (8) implies that $\rho_t \to 1$, and so from (6) it follows that there is $t_1 > t_0$ large enough so that $G(\rho_t z)/G(z)^{\rho_t} > (1 + \varepsilon)^{-1/2}$ for all[18] $z > t > t_1$. Therefore, for all $t > t_1$, we have

$$U \geq H(\gamma_t t)[G(\rho_t X_2) + G(\rho_t Y_1)]$$
$$\geq H(\gamma_t t)(1 + \varepsilon)^{-1/2}[G(X_2)^{\rho_t} + G(Y_1)^{\rho_t}]$$
$$\geq H(\gamma_t t)(1 + \varepsilon)^{-1/2} 2^{1 - \rho_t}[G(X_2) + G(Y_1)]^{\rho_t}.$$

---

[17] "$f(u, \varepsilon) = \varepsilon O(u^2)$ as $u \to 0$" means that there exists a constant $M < \infty$ such that $\overline{\lim}_{u \to 0} f(u, \varepsilon)/u^2 < \varepsilon M$ for every $\varepsilon \in (0, 1)$ [or, equivalently, for every $\varepsilon \in (0, 1)$ there exists $\underline{u} \equiv \underline{u}(\varepsilon) > 0$ such that $f(u, \varepsilon)/u^2 < \varepsilon M$ for all $u < \underline{u}$].

[18] Indeed, given $\varepsilon > 0$, let $z_0(\varepsilon)$ and $c_0(\varepsilon)$ be such that $G(cz)/G(z)^c > (1 + \varepsilon)^{-1/2}$ for all $z > z_0$ and $c \in (1, c_0)$; take $t_1 > z_0$ such that $\rho_t < c_0$ for all $t > t_1$.

Now $H(\gamma_t t)2^{1-\rho_t} \to 1$ by (9); therefore, there exists $t_2 \geq t_1$ such that $H(\gamma_t t) \times 2^{1-\rho_t} \geq (1+\varepsilon)^{-1/2}$ for all $t > t_2$, and so

$$U \geq (1+\varepsilon)^{-1}[G(X_2) + G(Y_1)]^{\rho_t}.$$

Let $u_2$ correspond to $t_2$; then $U < u < u_2$ implies that

(14) $$G(X_2) \leq ((1+\varepsilon)u)^{1/\rho_t} - G(Y_1) = G(t)^{1/\rho_t} - G(Y_1).$$

Equation (14) provides a *lower* bound on $X_2$, and so

$$\mathbf{P}(P_2 \in L^2, U < u | P_1) \leq \mathbf{P}(X_2 \text{ satisfies } (14) | Y_1) \leq G(t)^{1/\rho_t} - G(Y_1).$$

Integrating over $Y_1$ in $(t, \infty)$, we have

$$\mathbf{P}(P_1 \in L^1, P_2 \in L^2, U < u) \leq \int_t^\infty (G(t)^{1/\rho_t} - G(y_1)) \, dF(y_1)$$

$$= G(t)^{1+1/\rho_t} - \tfrac{1}{2} G(t)^2,$$

since $\int_t^\infty G(y_1) \, dF(y_1) = -\int_t^\infty G(y_1) \, dG(y_1) = -[G(y_1)^2/2]_t^\infty = G(t)^2/2$. Now $1/\rho_t = (1-\gamma_t)/((1+\gamma_t)(1+\gamma_t^2)) \geq 1 - 2\gamma_t$, and so $G(t)^{1+1/\rho_t} \leq G(t)^{2-2\gamma_t} \sim G(t)^2$ by (10), which implies that there is $t_3 \geq t_2$ such that $G(t)^{1+1/\rho_t} < (1+\varepsilon/2)G(t)^2$ for all $t > t_3$. This yields

$$\mathbf{P}(P_1 \in L^1, P_2 \in L^2, U < u) \leq \left(1 + \frac{\varepsilon}{2}\right) G(t)^2 - \frac{1}{2} G(t)^2$$

$$= \frac{1}{2}(1+\varepsilon)G(t)^2$$

$$= \frac{1}{2}(1+\varepsilon)^3 u^2$$

$$\leq \frac{1}{2}(1+7\varepsilon)u^2$$

for all $t > t_3$. Now let $u_3$ correspond to $t_3$. $\quad\square$

LEMMA 28. *As $u \to 0$,*

$$\mathbf{P}(P_1, P_2 \in L^1_{u,\varepsilon}, U < u) = \mathbf{P}(P_1, P_2 \in L^2_{u,\varepsilon}, U < u) = \varepsilon O(u^2).$$

PROOF. Let $P_1, P_2 \in L^1$ be such that $U < u < u_0$, where $u_0$ is given by Lemma 23, and let $t_0$ correspond to $u_0$. Then $Y_1 > Y_2 > t$ and $-\gamma_t t < X_1 < X_2 < \gamma_t t$, and also $C/B > t$ (by Lemma 23); therefore,

$$\frac{C}{A} = \frac{X_2 Y_1 - X_1 Y_2}{Y_1 - Y_2} \leq \frac{\gamma_t t Y_1 - (-\gamma_t t)Y_2}{Y_1 - Y_2} = \gamma_t t \frac{Y_1 + Y_2}{Y_1 - Y_2},$$

from which it follows by Lemma 26 with $\theta = \gamma_t$ that

$$U \geq H\left(\gamma_t \frac{C}{B}\right) G\left((1+\gamma_t)\frac{C}{A}\right) \geq H(\gamma_t t) G\left((1+\gamma_t)\gamma_t t \frac{Y_1 + Y_2}{Y_1 - Y_2}\right).$$

Let $t_1 > t_0$ be large enough so that $H(\gamma_t t) > 1/(1 + \varepsilon)$ for all $t > t_1$, and let $u_1$ correspond to $t_1$; then $U < u < u_1$ implies that

$$G\left((1 + \gamma_t)\gamma_t t \frac{Y_1 + Y_2}{Y_1 - Y_2}\right) < (1 + \varepsilon)u = G(t);$$

thus

$$(1 + \gamma_t)\gamma_t t \frac{Y_1 + Y_2}{Y_1 - Y_2} > t,$$

or

$$Y_1 < \rho_t Y_2,$$

where now $\rho_t := (1 + \gamma_t(1 + \gamma_t))/(1 - \gamma_t(1 + \gamma_t)) > 1$ and $\rho_t \to 1$. Therefore, for $P_2 \in L^1$,

$$\begin{aligned}
\mathbf{P}(P_1 \in L^1, U < u | P_2) &\le \mathbf{P}(Y_2 < Y_1 < \rho_t Y_2 | Y_2) \\
&= G(Y_2) - G(\rho_t Y_2) \\
&\le G(Y_2) - (1 - \varepsilon)G(Y_2)^{\rho_t},
\end{aligned}$$

the last inequality holding for all $t$ large enough, say $t > t_2 \ge t_1$, again by (6). Integrating over $Y_2$ in $(t, \infty)$ yields

$$\begin{aligned}
\mathbf{P}(P_1 \in L^1, P_2 \in L^1, U < u) &\le \int_t^\infty (G(y_2) - (1 - \varepsilon)G(y_2)^{\rho_t})\, dF(y_2) \\
&= \frac{1}{2}G(t)^2 - (1 - \varepsilon)\frac{1}{1 + \rho_t}G(t)^{1 + \rho_t}.
\end{aligned}$$

Now $1 + \rho_t \sim 2 + 2\gamma_t \to 2$, and so $G(t)^{1 + \rho_t}/(1 + \rho_t) \sim G(t)^2/2$ [recall (10)]. Therefore, there is $t_3 \ge t_2$ such that, for all $t > t_3$,

$$\begin{aligned}
\mathbf{P}(P_1 \in L^1, P_2 \in L^1, U < u) &\le \tfrac{1}{2}G(t)^2 - (1 - 2\varepsilon)\tfrac{1}{2}G(t)^2 \\
&= \varepsilon G(t)^2 = \varepsilon(1 + \varepsilon)^2 u^2 \le 4\varepsilon u^2.
\end{aligned}$$

The case where $P_1, P_2 \in L^2$ is the same (interchange $X$ and $Y$). $\quad\square$

LEMMA 29. *As* $u \to 0$,

$$\mathbf{P}(L_{u,\varepsilon}^3) = \varepsilon O(u).$$

PROOF. Define $L^4 := \{(x, y) : x + y > t\} \setminus (L^1 \cup L^2)$ (see Figure 5); then

$$\begin{aligned}
\mathbf{P}(L^4) &\le \mathbf{P}(x + y > t) - \mathbf{P}(|x| < \gamma_t t, y > (1 + \gamma_t)t) \\
&\quad - \mathbf{P}(|y| < \gamma_t t, x > (1 + \gamma_t)t).
\end{aligned}$$

Let $t_1$ be large enough so that $\mathbf{P}(x+y > t) \leq (2+\varepsilon)G(t)$ and also $H(\gamma_t t)G((1+\gamma_t)t) \geq (1-\varepsilon)G(t)$ for all $t > t_1$ [recall (7), (9) and (11)], and thus

$$\mathbf{P}(L^4) \leq (2+\varepsilon)G(t) - 2(1-\varepsilon)G(t) = 3\varepsilon G(t);$$

therefore,

$$\mathbf{P}(L^3) \leq \mathbf{P}(L^4) + \mathbf{P}(x \leq -\gamma_t t, y > t) + \mathbf{P}(x > t, y < -\gamma_t t)$$
$$\leq 3\varepsilon G(t) + 2F(-\gamma_t t)G(t) \leq 4\varepsilon G(t) = 4\varepsilon(1+\varepsilon)u$$

[note that $F(-\gamma_t t) \leq 1 - H(\gamma_t t) \leq \varepsilon$ for $t > t_1$]. $\quad\square$

COROLLARY 30. *As* $u \to 0$,

$$\mathbf{P}(P_1 \notin L^1_{u,\varepsilon} \text{ or } P_2 \notin L^2_{u,\varepsilon}, U < u) = \varepsilon O(u^2).$$

PROOF. For $u$ small enough,

$$\mathbf{P}(P_1 \notin L^1 \text{ or } P_2 \notin L^2, U < u)$$
$$= \mathbf{P}(P_1 \notin L \text{ or } P_2 \notin L, U < u)$$
$$\quad + \mathbf{P}(P_1 \in L^3, P_2 \in L, U < u) + \mathbf{P}(P_1 \in L, P_2 \in L^3, U < u)$$
$$\quad + \mathbf{P}(P_1, P_2 \in L^1, U < u) + \mathbf{P}(P_1, P_2 \in L^2, U < u)$$
$$\leq 0 + 2\mathbf{P}(L^3)\mathbf{P}(L) + 2\varepsilon O(u^2) = 2\varepsilon O(u)O(u) + 2\varepsilon O(u^2)$$

by Corollary 24 and Lemmas 29, 22 and 28. $\quad\square$

PROOF OF PROPOSITION 21. Adding up the estimates of Lemma 27 and Corollary 30 yields $\mathbf{P}(U < u) \leq (1/2)u^2 + \varepsilon O(u^2)$ as $u \to 0$. This holds for every $\varepsilon \in (0,1)$ and the left-hand side is independent of $\varepsilon$, and so[19] $\mathbf{P}(U < u) \leq (1/2)u^2 + o(u^2)$. $\quad\square$

5.3. *Second moment and Poisson approximation.* Recall Section 2.4 and Proposition 7 there.

PROPOSITION 31. *As* $n \to \infty$,

$$\mathbf{P}(\Gamma \cap \Gamma') = o(n^{-3}).$$

---

[19]Formally, there exists $M < \infty$ such that $\limsup_{u \to 0} \mathbf{P}(U < u)/u^2 \leq 1/2 + M\varepsilon$ for all $\varepsilon$, and so $\limsup_{u \to 0} \mathbf{P}(U < u)/u^2 \leq 1/2$.

PROOF. Proposition 32 below will show that $\mathbf{P}(W < u) = o(u^3)$. Thus, given $\varepsilon > 0$, there is $\delta > 0$ such that $\mathbf{P}(W < u) \leq \varepsilon u^3$ for all $u < \delta$; then, as in the proof of Theorem 19,

$$\int_0^1 (1-u)^{n-4} F_W(u)\, du \leq (1-\delta)^{n-4} + \varepsilon \int_0^1 (1-u)^{n-4} u^3\, du = \varepsilon O(n^{-4}).$$

Multiplying by $n-3$ and recalling that $\varepsilon > 0$ was arbitrary shows that indeed $\mathbf{P}(\Gamma \cap \Gamma') = o(n^{-3})$. $\square$

It remains to show that:

PROPOSITION 32. *As* $u \to 0$,
$$\mathbf{P}(W < u) = o(u^3).$$

PROOF. Fix $\varepsilon \in (0, 1)$. First, we have
$$\mathbf{P}(P_1 \in L^1, P_2 \in L^2, U < u, P_1' \in L^1, P_3' \in L^2, U' < u)$$
$$\leq \mathbf{P}(Y_1 > t, X_2 > t, Z_1 > t, X_3 > t) = G(t)^4 = O(u^4).$$
Next, for all $u$ small enough [i.e., $u < \underline{u}(\varepsilon)$],
$$\mathbf{P}(P_1 \notin L^1 \text{ or } P_2 \notin L^2, U < u, U' < u)$$
$$\leq \mathbf{P}(P_1 \notin L^1 \text{ or } P_2 \notin L^2, U < u, P_3' \in L) = \varepsilon O(u^2) O(u) = \varepsilon O(u^3)$$
by Corollary 30 and Lemma 22, and the fact that $P_3' = (X_3, Z_3)$ is independent of $P_1 = (X_1, Y_1)$ and $P_2 = (X_2, Y_2)$). Similarly, we have
$$\mathbf{P}(P_1' \notin L^1 \text{ or } P_3' \notin L^2, U' < u, U < u) = \varepsilon O(u^3).$$
Adding up the two terms yields $\mathbf{P}(W < u) \leq \varepsilon O(u^3)$ for every $\varepsilon \in (0, 1)$, or $\mathbf{P}(W < u) \leq o(u^3)$. $\square$

We can now prove Theorem 2.

PROOF OF THEOREM 2. Again, we apply Theorem 8 to $\sum_{i<j} D_{ij}$. We have
$$b_1 = O(n^2) O(n) \mathbf{P}(\Gamma)^2 = O(1/n),$$
$$b_2 = O(n^2) O(n) \mathbf{P}(\Gamma \cap \Gamma') = o(1),$$
by Theorem 19 and Proposition 31, and so $\|\mathcal{L}(S_2^{(n)}) - \text{Poisson}(\mu_n)\| \leq 2(b_1 + b_2) = o(1)$. Now $\text{Poisson}(\mu_n)$ converges to $\text{Poisson}(1/2)$ since $\mu_n \to 1/2$ by Theorem 19. $\square$

Recall (Proposition 4) that $S_1^{(n)}$, the number of one-point (pure) ESS, converges in distribution to Poisson(1) as $n \to \infty$. While $S_1^{(n)}$ and $S_2^{(n)}$ are not independent, we will now show that, nevertheless, their sum converges to Poisson(3/2).

THEOREM 33. *Put* $S_{\leq 2}^{(n)} := S_1^{(n)} + S_2^{(n)}$. *If* $F \in \mathcal{SE}$, *then*

$$\mathcal{L}(S_{\leq 2}^{(n)}) \to \text{Poisson}(3/2) \qquad \text{as } n \to \infty.$$

PROOF. We apply again Theorem 8, this time to $\sum_i C_i + \sum_{i<j} D_{ij}$. Let $b_1'$ and $b_2'$ correspond to $S_2^{(n)} = \sum_{i<j} D_{ij}$; in the proof of Theorem 2 above we showed that $b_1' = O(1/n)$ and $b_2' = o(1)$. The additional dependencies now are between a term $C_i$ and a term $D_{ij}$, with the same $i$ and $j \neq i$. However, we have $\mathbf{E}(C_i D_{ij}) = \mathbf{P}(C_i = D_{ij} = 1) = 0$, since $C_i = 1$ implies that $R_{ii} > R_{ij}$, whereas $D_{ij} = 1$ implies that $R_{ii} < R_{ij}$ (see Lemma 3). Thus $b_2 = b_2' = o(1)$, and

$$b_1 = b_1' + \sum_i \mathbf{P}(C_i = 1)^2 + 2 \sum_i \sum_{j \neq i} \mathbf{P}(C_i = 1)\mathbf{P}(D_{ij} = 1)$$

$$= O(1/n + n(1/n)^2 + n^2(1/n)(1/n^2)) = O(1/n).$$

Theorem 8 yields $\|\mathcal{L}(S_{\leq 2}^{(n)}) - \text{Poisson}(1 + \mu_n)\| \leq 2(b_1 + b_2) = o(1)$; and we have $1 + \mu_n \to 3/2$ by Theorem 19. □

COROLLARY 34. *If* $F \in \mathcal{SE}$, *then the probability that there is an ESS with support of size* $\leq 2$ *converges to* $1 - e^{-3/2} \simeq 0.78$ *as* $n \to \infty$.

**6. Discussion.** We conclude with a discussion of some of the related literature, together with a number of comments, conjectures and open problems.

(a) *Vertices and equilibria.* The connection between Nash equilibria and vertices of random polytopes was used by Bárány, Vempala and Vetta [3] to find Nash equilibria in random games. Concerning ESS, we emphasize again that the number of vertices of a random polygon and the number of two-point ESS of a random game have different distributions; only their expectations are related (by a factor of 8; see Proposition 9).

(b) *The class $\mathcal{EF}$.* The class of distributions with "Exponential and Faster decreasing tails" for which Theorem 1 holds can clearly be taken to be larger than that of Section 4. Indeed, since Theorem 13 holds for *any* distribution, we can include in $\mathcal{EF}$ any $F$ such that $\mu_n \to \infty$. Take, for example, those distributions in Fisher [8] for which the limit shape of the convex hull is a *strictly* convex set; this implies that the number of vertices, and thus $\mu_n$, must go to infinity. By Theorem 1 there, this includes distributions where, for some $\alpha > 1$, the tail probability $G = 1 - F$ satisfies $G^{-1}(1/t^c) \sim c^{1/\alpha} G^{-1}(1/t)$ as $t \to \infty$ for each $c \in (0, 1)$.

(c) *The probability of having an ESS.* For distributions in $\mathcal{EF}$, Theorem 1(iii) implies that the probability that there is an ESS converges to 1 as $n$ increases. For distributions in $\mathcal{SE}$, however, it is still unknown what the limit of this probability is. Some preliminary informal analysis suggests to us the following conjectures: if $F \in \mathcal{SE}$, then, as $n \to \infty$:

- $S_3^{(n)} \to \text{Poisson}(1/3)$;
- $S_{\leq 3}^{(n)} := S_1^{(n)} + S_2^{(n)} + S_3^{(n)} \to \text{Poisson}(1 + 1/2 + 1/3) = \text{Poisson}(11/6)$;
- $S_\ell^{(n)} \to 0$ for all $\ell \geq 4$;
- $\sum_{\ell=1}^n S_\ell^{(n)} \to \text{Poisson}(11/6)$;
- $\mathbf{P}(\text{there is an ESS}) \to 1 - e^{-11/6} \simeq 0.84 < 1$.

[The geometric objects corresponding to $S_\ell^{(n)}$ are now the $(\ell-1)$-dimensional faces of the convex hull of $n$ random points in $\mathbb{R}^\ell$.]

(d) *Threshold phenomenon.* Our distributions exhibit a "threshold" phenomenon: either $\mu_n \to \infty$ or $\mu_n \to 1/2$. However, we believe that one may construct distributions for which the sequence $\mu_n$ has other limit points, or even oscillates wildly as $n$ increases. Indeed, for each $n$, the number of vertices, and thus $\mu_n$, depends on the distribution $F$ only through a certain interval of its tail [in a neighborhood of $G^{-1}(1/n)$]. Therefore, one should be able to "glue" various tails (of the $\mathcal{EF}$ or $\mathcal{SE}$ types) and get different limit points. See Devroye [7] for such oscillations in the case of radially symmetric distributions.

(e) *Other distributions.* It would be interesting to study additional classes of distributions. For example, bounded-support distributions whose tail $G$ is not log-concave are not included in $\mathcal{EF}$; we conjecture that $\mu_n \to \infty$ in this case, though perhaps the convergence is at a slower rate than the $\log n$ of the uniform distribution. Another question arises when the distributions of the $X$-coordinates and of the $Y$-coordinates differ (but, say, they are both in $\mathcal{EF}$ or both in $\mathcal{SE}$); this concerns also the number of vertices when the distribution is not symmetric (consider the different orthants).

**Acknowledgments.** The authors thank Imre Bárány, Josef Hofbauer, Nati Linial, Noam Nisan, Bernhard von Stengel and two anonymous referees for helpful comments.

S. HART
INSTITUTE OF MATHEMATICS
DEPARTMENT OF ECONOMICS
AND
CENTER FOR THE STUDY OF RATIONALITY
HEBREW UNIVERSITY OF JERUSALEM
ISRAEL
E-MAIL: hart@huji.ac.il
URL: http://www.ma.huji.ac.il/hart

Y. RINOTT
DEPARTMENT OF STATISTICS
AND
CENTER FOR THE STUDY OF RATIONALITY
HEBREW UNIVERSITY OF JERUSALEM
ISRAEL
E-MAIL: rinott@mscc.huji.ac.il

B. WEISS
INSTITUTE OF MATHEMATICS
AND
CENTER FOR THE STUDY OF RATIONALITY
HEBREW UNIVERSITY OF JERUSALEM
ISRAEL
E-MAIL: weiss@math.huji.ac.il